\documentclass[12pt]{amsart}
\usepackage[T1]{fontenc}
\usepackage[utf8]{inputenc}
\usepackage{lmodern}
\usepackage{amsmath,amssymb,amsthm,mathtools}
\usepackage{geometry}
\usepackage{microtype}
\usepackage[hidelinks]{hyperref}
\numberwithin{equation}{section}
\newtheorem{theorem}{Theorem}[section]
\newtheorem{proposition}{Proposition}[section]
\newtheorem{lemma}[proposition]{Lemma}
\newtheorem{corollary}[proposition]{Corollary}
\newtheorem{citedtheorem}[proposition]{Theorem}

\theoremstyle{remark}
\newtheorem{remark}[proposition]{Remark}
\newcommand{\R}{\mathbb R}

\newcommand{\Z}{\mathbb Z}
\newcommand{\Sn}{\mathbb S}
\newcommand{\Vol}{\operatorname{Vol}}
\newcommand{\ind}{\operatorname{ind}}
\newcommand{\tr}{\operatorname{tr}}
\newcommand{\rank}{\operatorname{rank}}
\newcommand{\Cat}{\operatorname{Cat}}
\newcommand{\dd}{\,\mathrm d}
\newcommand{\ip}[2]{\langle #1,#2\rangle}
\newcommand{\norm}[1]{\lvert #1\rvert}
\title[Sharp signature-weighted total-curvature bounds]{Sharp signature-weighted total-curvature bounds for complete minimal submanifolds}
\author{Jiuru Zhou}
\address{School of Mathematics, Yangzhou University,
	Yangzhou, 225002, China}
\email{zhoujiuru@yzu.edu.cn}
\thanks{Partially supported by
	NSFC Grants 12171417.}
\begin{document}
\maketitle

\begin{abstract}
Let $i:M^n\to\R^{n+p}$ be a complete connected minimal immersion,
with $n\ge3$. We get a sharp bound of $\int_M |A|^n\dd\mu$ from below
in terms of the number of ends and the mod-two Betti numbers of $M$. As an application, if
$\int_M |A|^n\dd\mu<\gamma_n$,
then $i(M)$ is an affine $n$-plane.
Moreover, if $\int_M |A|^n\dd\mu=\gamma_n$, then $i(M)$ is an $n$-catenoid
in an affine $(n+1)$-plane, where $\gamma_n=\frac{[n(n-1)]^{n/2}}{n-1}\Vol(\Sn^n)$.
\end{abstract}

\section{Introduction}\label{sec:intro}
Let $i:M^n\to\R^{n+p}$ be a complete connected minimal immersion.
By the Gauss equation the scalar curvature of $M$ is $-|A|^2$. Then 
\begin{equation}\label{eq:energy}
\int_M |A|^n\dd\mu
\end{equation}
is called the total scalar curvature. 
Anderson \cite[Corollary 2.3]{Anderson} showed that there is a constant $\epsilon_0(n,p)$ such that if the total scalar curvature $\int_M|A|^n\dd\mu\leq\epsilon_0(n,p)$, then $i(M)$ is an affine $n$-plane. For $n\geq 4,~q=1$, Ni \cite[Theorem A]{Ni} proved that there is a constant $C_1(n)$ such that if $(\int_M|A|^n\dd\mu)^{\frac{1}{n}}<C_1(n)$, then $M$ must be totally geodesic, and also get small $(\int_M|A|^n\dd\mu)^{\frac{1}{n}}$ controls the number of ends. Later, Seo obtained a positive gap $C(n)$ of the total scalar curvature in arbitrary codimension for $n\ge3$
\cite[Corollary~1.2]{Seo}. These results suggest asking how much
curvature each additional end forces, and which immersions attain
the resulting bound.

For a compact immersion, the famous Chern-Lashof theory related total absolute curvature
to the average number of critical points of a height function
\cite{ChernLashof}. Later, Shiohama-Xu gives a separate integral
for each signature of the shape operator \cite{ShiohamaXu1997,ShiohamaXu2011}.
Vlachos used this relation to bound intrinsic curvature integrals by
topology for compact Euclidean submanifolds \cite{Vlachos2013}.
For complete noncompact immersions, van Gemmeren developed a relative
version of the Morse inequalities \cite{vanGemmeren1998}. Here we want to apply Chern-Lashof type theory to minimal submanifolds.

Let's denote $e=e(M)$ the number of ends and
$b_k=b_k(M)=\dim H_k(M;\Z_2)$. Set
\begin{equation}\label{eq:gamma}
 \gamma_n=\frac{[n(n-1)]^{n/2}}{n-1}\Vol(\Sn^n),\qquad
 c_{n,k}=n^{n/2}\left(\frac{n-k}{k}\right)^{n/2-k}
 \quad (1\le k\le n-1).
\end{equation}

\begin{theorem}\label{main1}
Let $i:M^n\to\R^{n+p}$ be a complete connected minimal immersion without
boundary, with $n\ge3$, $p\ge1$ and $\int_M |A|^n\dd\mu<\infty$. Then
\begin{equation}\label{eq:main-topology}
\begin{aligned}
 \int_M |A|^n\dd\mu\ge {}&\gamma_n(b_1+e-1)\\
 &+\Vol(\Sn^n)\left[\sum_{2\le k<n/2}c_{n,k}b_k+
 \begin{cases}
  \tfrac12 n^{n/2}b_{n/2},& n\text{ even},\\
  0,& n\text{ odd}.
 \end{cases}\right],
\end{aligned}
\end{equation}
with optimal coefficient for every fixed $p$.
\end{theorem}

\begin{remark}
	(1). In particular, we have
	\begin{equation}\label{eq:topology-ends}
		\int_M |A|^n\dd\mu\ge\gamma_n(b_1+e-1),
	\end{equation}
	so
	\begin{equation}\label{eq:ends}
		\int_M |A|^n\dd\mu\ge\gamma_n(e-1),\qquad e\le1+\frac{\int_M |A|^n\dd\mu}{\gamma_n}.
	\end{equation}

(2). The ordinary catenoid has $b_1=0$,
$e=2$ and $\int_M |A|^n\dd\mu=\gamma_n$, so it realizes equality in
\eqref{eq:topology-ends}. Geometric optimality of the individual
higher Betti-number coefficients remains a separate question.
\end{remark}

As an application, we obtain the following sharp gap and rigidity results for the total scalar curvature,
\begin{theorem}\label{main2}
Let $i:M^n\to\R^{n+p}$ be a complete connected minimal immersion without
boundary, with $n\ge3$ and $p\ge1$. If
\begin{equation}\label{eq:sharp-gap}
 \int_M |A|^n\dd\mu<\gamma_n,
\end{equation}
then $i(M)$ is an affine $n$-plane.
Moreover, if $\int_M |A|^n\dd\mu=\gamma_n$, then $i(M)$ is an $n$-catenoid
in an affine $(n+1)$-plane.
\end{theorem}

The proof begins with a generic linear height function. On each
graphical end a sufficiently low sublevel is contractible. Thus the low sublevel $L$
has $e$ components, and the exact sequence of $(M,L)$ gives
$\dim H_1(M,L;\Z_2)=b_1+e-1$. Relative Morse theory bounds this rank
by the number of index-one critical points. The normal Gauss map turns
the average of this count into an integral of the absolute determinant
of the shape operator. Applying the trace-free determinant inequality
on each index stratum, then averaging over the normal sphere, gives
Theorem \ref{main1}. At energy $\gamma_n$ the inequality in Theorem \ref{main1} and one-end rigidity force two ends. Then
equality in the normal average forces a one-dimensional first normal
space wherever $A\ne0$. The Codazzi equation reduces the immersion to
a hypersurface, where Schoen's two-end theorem applies.

Section~\ref{sec:infinity} collects the end structure, the homology of
a compact core, the normal Gauss-map formulas and the two-end
classification theorems. Section~\ref{sec:morse}
computes the relative Morse bounds on the ends. Section~\ref{sec:algebra}
proves the two algebraic inequalities and Theorem~\ref{main1}. In
Section~\ref{sec:rigidity} we compute the curvature of the ordinary
catenoid, prove Theorem~\ref{main2}.

\section{Preliminaries}\label{sec:infinity}
In this section, we list some results, most of which are well-known.
Let $N=n+p$ and let $g$ be the metric induced by $i:M^n\to\R^N$.
Throughout the paper, when performing ambient Euclidean operations, we
write $x$ for $i(x)\in\R^N$ and identify $X\in T_xM$ with
$di_x(X)\in\R^N$. If $D$ denotes the Euclidean
connection and $\nabla$ the
Levi--Civita connection of $g$, then
\[
 A(X,Y)=D_XY-\nabla_XY,\qquad
 \operatorname{Scal}_M=|\tr_g A|^2-|A|^2.
\]
For a minimal immersion $\tr_g A=0$ and hence
$\operatorname{Scal}_M=-|A|^2$. 

\subsection{Finite total curvature and graphical ends}\label{sec:anderson}

\begin{proposition}[Finite-curvature end structure]\label{prop:end-structure}
Suppose $i:M^n\to\R^N$ is a complete connected minimal immersion
without boundary, $n\ge3$, and $\int_M|A|^n\dd\mu<\infty$.
Fix $o\in M$, translate the ambient coordinates so that $i(o)=0$, and set
$s(x)=d_M(o,x)$ and $\rho(x)=|x|$. Then
\[
 s|A|\longrightarrow0,\qquad \frac{\rho}{s}\longrightarrow1,
 \qquad |\nabla\rho|\longrightarrow1
\]
uniformly as $s\to\infty$. The immersion is proper. For all sufficiently
large $R_0$, $\{\rho\ge R_0\}$ is diffeomorphic to
$\rho^{-1}(R_0)\times[R_0,\infty)$ with second coordinate $\rho$.
There are finitely many ends, and on each end the tangent planes
converge uniformly to a single unoriented linear $n$-plane.
\end{proposition}

\begin{proof}
	By Anderson \cite[Proposition~2.2]{Anderson},
	$s|A|\to0$ uniformly as $s\to\infty$.
	Fix $0<\varepsilon<1$ and choose $T>0$ such that
	$s|A|\le\varepsilon$ on $\{s\ge T\}$.
	For $x\in M$ with $s=s(x)>T$, let
	$\gamma:[0,s]\to M$ be a unit-speed minimizing geodesic from
	$o$ to $x$, and put $X(t)=i(\gamma(t))$.
	Since $i(o)=0$, we have
	$|X'(t)|=1$, $|X(t)|\le t$, and, by the Gauss formula,
	$X''(t)=A(\gamma'(t),\gamma'(t))$.
	Consequently, for $T\le t\le s$,
	\[
	\frac{d}{dt}\langle X(t),X'(t)\rangle
	=1+\langle X(t),X''(t)\rangle
	\ge1-t|A(\gamma(t))|\ge1-\varepsilon.
	\]
	As $\langle X(T),X'(T)\rangle\ge-T$, integration yields
	\[
	(1-\varepsilon)s-(2-\varepsilon)T
	\le\langle X(s),X'(s)\rangle\le\rho(x)\le s.
	\]
	It follows, first letting $s\to\infty$ and then
	$\varepsilon\downarrow0$, that $\rho/s\to1$ uniformly.
	Moreover, since $X'(s)$ is a unit tangent vector and
	$\nabla\rho=(X(s)/\rho)^\top$, for sufficiently large $s$,
	\[
	1\ge|\nabla\rho|
	\ge\frac{\langle X(s),X'(s)\rangle}{\rho(x)}
	\ge1-\varepsilon-\frac{(2-\varepsilon)T}{s}.
	\]
	Thus $|\nabla\rho|\to1$ uniformly as well. The convergence
	$\rho/s\to1$ implies that $\rho\ge s/2$ outside a compact
	intrinsic ball; by Hopf--Rinow, inverse images of bounded
	Euclidean closed balls are compact. Hence $i$ is proper.
	
	Choose $R_0$ large enough that $|\nabla\rho|\ge1/2$
	on $\{\rho\ge R_0\}$. Properness makes
	$\Sigma_{R_0}=\rho^{-1}(R_0)$ compact, and the normalized
	radial field
	\[
	V=\frac{\nabla\rho}{|\nabla\rho|^2}
	\]
	is smooth on $\{\rho\ge R_0\}$, satisfies $|V|\le2$, and
    $V(\rho)=1$. Since $M$ is complete, its flow
	$\Phi_t$ cannot escape in finite time. It therefore exists
	forward for all $t\ge0$, and backward until $\rho=R_0$.
	As $\rho(\Phi_t(x))=\rho(x)+t$, the map
	\[
	\Sigma_{R_0}\times[R_0,\infty)\longrightarrow
	\{\rho\ge R_0\},\qquad
	(y,r)\longmapsto\Phi_{r-R_0}(y)
	\]
	is a diffeomorphism whose second coordinate is $\rho$;
	compare \cite[Lemma~2.4 and Corollary~2.5]{Anderson}.
	The compact manifold $\Sigma_{R_0}$ has finitely many
	connected components, each determining one end. Since
	$\Delta_M|i|^2=2n$, $M$ cannot be compact, so the number
	of ends is positive and finite.
	
	Finally, Anderson \cite[Theorems~3.1--3.2]{Anderson}
	shows that the tangent-plane Gauss map extends continuously
	to a compactification obtained by adding one point for each
	end. Its value at that point is a unique unoriented limiting
	$n$-plane. 
	By continuity of this extension, the tangent planes
	converge uniformly to the limiting plane on each end.
	More precisely, if $P_j$ denotes the limiting $n$-plane
	of the end $E_j$, then
	\[
	\lim_{R\to\infty}
	\sup_{\substack{x\in E_j\\ s(x)\ge R}}
	d_{\mathrm{Gr}}\bigl(di_x(T_xM),P_j\bigr)=0,
	\]
	where $d_{\mathrm{Gr}}$ denotes a metric inducing the
	standard topology on the Grassmannian
	$\mathrm{Gr}(n,N)$.
\end{proof}

\begin{proposition}[Graphical ends]\label{prop:graph-ends}
Under the hypotheses of Proposition~\ref{prop:end-structure}, $M$
retracts by deformation onto a connected compact core $K$. Each
component of $M\setminus K$ is mapped by $i$ diffeomorphically onto
a graph over the exterior of a ball in its limiting plane $P$.
In coordinates in $P\oplus P^\perp$, the graph has the form
\begin{equation}\label{eq:graph}
 x\longmapsto(x,u(x)),\qquad |x|>R,\qquad
 |Du(x)|\to0,\quad u(x)=o(|x|).
\end{equation}
Both estimates hold uniformly on the end.
\end{proposition}

\begin{proof}	
	Fix an end $E$ given by
	Proposition~\ref{prop:end-structure}, and let $P$ be its
	limiting $n$-plane. By the uniform convergence of the tangent
	planes to $P$, after increasing $R_0$ if necessary, the
	orthogonal projection
	\[
	F:=\pi_P\circ i:E\longrightarrow P
	\]
	is a local diffeomorphism.
	
	Since $\rho(x)=|i(x)|$, the intrinsic gradient
	of $\rho$ is the tangential projection of the ambient
	radial unit vector:
	\[
	\nabla\rho
	=\left(\frac{i(x)}{\rho(x)}\right)^\top.
	\]
	Since $|i(x)/\rho(x)|=1$, orthogonal decomposition
	into tangent and normal components gives
	\[
	\operatorname{dist}\left(
	\frac{i(x)}{\rho(x)},T_xM\right)^2
	=
	\left|\left(\frac{i(x)}{\rho(x)}\right)^\perp\right|^2
	=1-|\nabla\rho(x)|^2
	\longrightarrow0
	\]
	uniformly as $s(x)\to\infty$.
	
	Since $T_xM\to P$ uniformly, it follows that
	\[
	\frac{|\pi_{P^\perp}i(x)|}{\rho(x)}\longrightarrow0,
	\qquad
	\frac{|F(x)|}{\rho(x)}\longrightarrow1.
	\]
	In particular, $|F(x)|\ge\rho(x)/2$ for sufficiently large $x$.
	
	Choose $R$ sufficiently large such that
	$|F|<R$ on $\partial E$ and $|F(x)|\ge \frac{1}{2}\rho(x)~(x\in E,~\rho(x)>2R)$. Define
	$U:=E\cap\{|F|>R\}$.
	The restriction
	\[
	F:U\longrightarrow P\setminus\overline B_R
	\]
	is a local diffeomorphism. For any compact set $C\subset P\textbackslash\bar{B}_R$, we next show that $F^{-1}(C)$ is compact. Since $C$ is compact, there exists a constant $B>R$ such that $|y|\leq B~(y\in C)$. If $x\in F^{-1}(C)$, $|F(x)|\leq B$. In addition, $|F(x)|\ge \frac{1}{2}\rho(x)$ for large $x$. Hence, 
	\[
	 \rho(x)\leq 2|F(x)|\leq 2B.
	\]
	There exists a constant $B'>0$ such that 
	\[
	 F^{-1}(C)\subset\{x\in M:\rho(x)\leq B'\}.
	\]
	Since the immersion $i$ is proper, the right hand side is compact. Since $|F|<R$ on $\partial E$ and $|y|>R$ for $y\in C$, $F^{-1}(C)\cap \partial E=\emptyset$. Hence, $F$ is proper. Then $F:~U\to P\setminus\overline B_R$ is a covering.
	
	Suppose the connected components of $U$ are $\{U_\alpha\}_{\alpha\in I}$. The each 
	$F|_{U_\alpha}:~U_\alpha\to P\setminus\overline B_R$ is proper local diffeomorphism and $F(U_\alpha)$ is open and closed. Hence, $F(U_\alpha)=P\textbackslash \bar{B}_R$, since $P\textbackslash \bar{B}_R$ is connected. Moreover, $U_\alpha$ is unbounded. For any fixed $R>0$, there exists $x_\alpha\in U$ with $\rho(x_\alpha)>2R$. Hence every
	component $U_\alpha$ meets $E\cap\{\rho>2R\}$, which is connected by
	Proposition~\ref{prop:end-structure} and lies in $U$.
	Thus $U$ is connected. Since $n\ge3$, $P\textbackslash \bar{B}_R$ is simply
	connected.
	Thus, the covering has one sheet and $F$ is a
	diffeomorphism.
	Writing its inverse in $P\oplus P^\perp$ and define $u(x)=\pi_{P^\perp}(i(F^{-1}(x)))$, gives
	\[
	i\circ F^{-1}(x)=(x,u(x)),\qquad |x|>R.
	\]
	The convergence $T_xM\to P$ implies $|Du(x)|\to0$
	uniformly. Moreover, for $r>R$ and $\theta\in\Sn^{n-1}$,
	\[
	|u(r\theta)|
	\le \sup_{|\xi|=R}|u(\xi)|
	+\int_R^r |Du(t\theta)|\,dt.
	\]
	Dividing by $r$ and using the uniform decay of $Du$
	gives $u(x)=o(|x|)$ uniformly as $|x|\to\infty$.

	Finally, we construct the compact core and the deformation
	retraction. Perform the preceding construction on all $e$
	ends $E_1,\ldots,E_e$. Let
	$F_j=\pi_{P_j}\circ i$ denote the orthogonal projection
	on the $j$th end. Since there are only finitely many ends,
	we may choose a common sufficiently large radius $R$
	such that the graphical parametrizations are defined
	on a neighborhood of $\{|x|\ge R\}$ and the corresponding
	tails are pairwise disjoint.
	
	For each $j$, define the closed graphical tail
	\[
	T_j=\{F_j^{-1}(r\theta):
	r\ge R,\ \theta\in\mathbb S^{n-1}\}.
	\]
	Since $F_j$ is a diffeomorphism on the graphical end,
	we have
	\[
	T_j\cong\mathbb S^{n-1}\times[R,\infty),
	\]
	whose boundary is the sphere
	\[
	\Sigma_j
	=\{F_j^{-1}(R\theta):
	\theta\in\mathbb S^{n-1}\}
	\cong\mathbb S^{n-1}.
	\]
	
	Remove the open tails and set
	\[
	K=M\setminus\bigcup_{j=1}^e
	\{F_j^{-1}(r\theta):r>R\}.
	\]
	Then $K$ is a smooth manifold with boundary
	\[
	\partial K=\bigsqcup_{j=1}^e\Sigma_j
	\cong\bigsqcup_{j=1}^e\mathbb S^{n-1}.
	\]
	Moreover, since $|F_j(q)|/\rho(q)\to1$ uniformly
	on each end, every point of $M$ with sufficiently
	large $\rho$ belongs to one of the removed tails.
	Thus $K$ is contained in $\{\rho\le B\}$ for some
	$B>0$. As $K$ is closed in $M$ and the immersion
	$i$ is proper, $K$ is compact.
	
	We now construct a strong deformation retraction
	of $M$ onto $K$. On each closed tail $T_j$, write
	$q=F_j^{-1}(r\theta)$, where $r\ge R$, and define
	\[
	H_t(q)
	=F_j^{-1}\bigl(((1-t)r+tR)\theta\bigr),
	\qquad 0\le t\le1.
	\]
	On $K$, set $H_t(q)=q$ for all $t\in[0,1]$.
	These definitions agree on each boundary sphere
	$\Sigma_j$, since
	\[
	H_t(F_j^{-1}(R\theta))
	=F_j^{-1}(R\theta).
	\]
	Hence they define a continuous map
	$H:M\times[0,1]\to M$ satisfying
	\[
	H_0=\operatorname{id}_M,\qquad
	H_1(M)=K,\qquad
	H_t|_K=\operatorname{id}_K.
	\]
	Therefore $H$ is a strong deformation retraction
	of $M$ onto $K$. In particular, $K=H_1(M)$ is
	connected because $M$ is connected.
\end{proof}

\begin{citedtheorem}[Anderson {\cite[Theorem 5.2]{Anderson}}]\label{lem:one-end}
Under the hypotheses of Proposition~\ref{prop:end-structure}, if $e(M)=1$, the immersion parametrizes an affine $n$-plane.
\end{citedtheorem}

\subsection{Homology of a compact core}

\begin{lemma}[Homology of the compact core]\label{rem:betti-symmetry}
Under the hypotheses of Theorem \ref{main1}, $M$ deformation retracts onto a
connected compact core $K$ with $e$ boundary spheres. With coefficients
in $\Z_2$,
\begin{equation}\label{eq:betti-end-duality}
 b_{n-1}(M)=b_1(M)+e-1,\qquad
 b_{n-k}(M)=b_k(M)\quad(2\le k\le n-2).
\end{equation}
\end{lemma}
\begin{proof}
Cut the graphs of Proposition~\ref{prop:graph-ends} at a large
coordinate sphere and retract each remaining exterior part onto its
boundary. Cap the resulting core $K$ with $e$ balls to form a closed
connected manifold $\widehat M$. By excision, $H_j(\widehat M,K)$
vanishes except for $j=n$, where it is $(\Z_2)^e$. The fundamental class
of $\widehat M$ maps to $(1,\ldots,1)$, one relative fundamental class
in each cap. The exact sequence of the pair therefore gives
\[
 b_{n-1}(K)=b_{n-1}(\widehat M)+e-1,\qquad
 b_j(K)=b_j(\widehat M)\quad(1\le j\le n-2).
\]
Since $n\ge3$, the second identity applies to $j=1$.
Poincar\'e duality on $\widehat M$ now gives \eqref{eq:betti-end-duality}.
\end{proof}

\subsection{Height functions and the normal Gauss map}\label{sec:height-geometry}
Following the Chern-Lashof \cite{ChernLashof} approach, for $a\in\Sn^{N-1}$ consider the linear height function $h_a(x)=\langle i(x),a\rangle$ and set
$M_a^t=\{h_a\le t\}$. For a unit normal $\nu$ at $x$ let $S_\nu$
be the shape operator defined by $\ip{S_\nu X}{Y}=\ip{A(X,Y)}{\nu}$.
By Gauss formula,
\begin{equation}\label{eq:hessian-height}
\nabla h_a=a^\top,\qquad
\nabla^2h_a(X,Y)=\ip{A(X,Y)}{a}.
\end{equation}
Thus $x$ is critical precisely when $a$ is normal at $x$, and the
Hessian there is represented by $S_a$. Since $\tr S_a=0$, every
nondegenerate critical point has index between $1$ and $n-1$.

Let $G:N^1M\to\Sn^{N-1}$ be the map $G(x,\nu)=\nu$ on the unit
normal bundle. The normal connection identifies its tangent space at
$(x,\nu)$ with $T_xM\oplus(N_xM\cap\nu^\perp)$. Give the two
summands their induced metrics and make them orthogonal. The projection
onto $M$ is then a Riemannian submersion with round fibers, so
\begin{equation}\label{eq:normal-fiber-integral}
\int_{N^1M} F\dd V
 =\int_M\int_{S(N_xM)}F(x,\nu)\dd\nu\dd\mu(x)
\end{equation}
for every nonnegative measurable $F$. Here $\dd\nu$ is round measure
on the unit normal sphere, with counting measure on $\Sn^0$ when $p=1$.
In this splitting the Weingarten formula gives
\begin{equation}\label{eq:gauss-jacobian}
 dG_{(x,\nu)}(X,\eta)=-S_\nu X+\eta,\qquad
 |\det dG_{(x,\nu)}|=|\det S_\nu|.
\end{equation}
The fiber over
$a$ consists of the pairs $(x,a)$ with $dh_a(x)=0$. Thus regular values
of $G$ give Morse heights by \eqref{eq:hessian-height}, in every
codimension and without any orientation assumption.

\subsection{Regularity at infinity and two-end rigidity}\label{sec:schoen-expansion}

\begin{proposition}[Hypersurface regularity at infinity]\label{prop:hypersurface-regularity}
Under the hypotheses of Proposition~\ref{prop:end-structure},
suppose that $i(M)$ lies in an affine $(n+1)$-plane
$V$. As a hypersurface of $V$, the immersion is regular at infinity.
In coordinates over the limiting plane of each end the scalar graph
function satisfies
\begin{equation}\label{eq:schoen-expansion}
 u(x)=b+c|x|^{2-n}+\sum_{j=1}^n d_jx_j|x|^{-n}+O(|x|^{-n}),
 \qquad |Du|=O(|x|^{1-n}),\quad |D^2u|=O(|x|^{-n}).
\end{equation}
\end{proposition}
\begin{proof}
Identify $V$ with $\R^{n+1}$. Proposition~\ref{prop:graph-ends}
gives finitely many scalar exterior graphs of bounded slope.
Schoen's criterion and the estimates in its proof then give
\eqref{eq:schoen-expansion} and regularity at infinity
\cite[Proposition~3, pp.~802--804]{Schoen1983}.
\end{proof}

\begin{citedtheorem}[Schoen {\cite[Theorem~3]{Schoen1983}}]\label{thm:schoen}
Let $n\ge3$. A complete minimal hypersurface immersion in $\R^{n+1}$ that is
regular at infinity and has two ends is a catenoid or a pair of
affine hyperplanes.
\end{citedtheorem}

\section{Relative Morse theory at infinity}\label{sec:morse}
Let $i:M^n\to\R^{n+p}$ satisfy the hypotheses of Theorem~\ref{main1} and
write $e=e(M)$. We first determine the topology of a low sublevel
of a generic height function. This gives a lower bound for the number of
critical points of each index. 

\subsection{Sublevels on the ends}
\begin{lemma}\label{lem:generic-height}
For almost every $a\in\Sn^{N-1}$, $h_a$ is Morse with finitely many
critical points and $|\nabla h_a|\ge c_a>0$ outside a compact set.
For all sufficiently negative $t$, $M_a^t$ has $e$ components, each
diffeomorphic to a closed half-space. For all sufficiently positive
$T$, $M_a^T$ is connected.
\end{lemma}
\begin{proof}
	By Proposition \ref{prop:graph-ends}, there exists a compact core $K\subset M$, such that $M\textbackslash K=E_1\cup\cdots\cup E_e$. For each end $E_j$, let $P_j$ be its limiting tangent plane and set $F_j=\pi_{P_j}\circ i|_{E_j}$. In these coordinates, $i\circ F_j^{-1}(x)=(x,u_j(x))$, $x\in P_j\setminus\overline B_{R_0}$. The graph functions satisfy 
	\[
	 |Du_j(x)|\to 0,~\frac{|u_j(x)|}{|x|}\to 0~(|x|\to\infty).
	\]
	Exclude the critical values of
	$G$ and the spheres $S(P_j^\perp)$ from $\Sn^{N-1}$. The first set
	has measure zero by Sard's theorem and each of the spheres has
	dimension $p-1<N-1$, also has measure zero.
	Fix a direction $a$ outside these sets.
	
	On the $j$th graphical end, write $a=a_{P_j}+a_{P_j^\perp}$.
	The excluded directions ensure that
	$q_j=\norm{a_{P_j}}>0$. Rotate coordinates in $P_j$ so that
	$a_{P_j}=q_je_1$. In these graph coordinates, the height function satisfies
	\[
	h_a(F_j^{-1}(x))=q_jx_1+v_j(x),
	\qquad
	v_j(x)=\text{constant}+\ip{u_j(x)}{a_{P_j^\perp}}.
	\]
	The vector $a_{P_j^\perp}$ is fixed, so
	\[
	\norm{Dv_j(x)}
	\le\norm{a_{P_j^\perp}}\norm{Du_j(x)}
	\le\norm{Du_j(x)}\longrightarrow0,
	\qquad v_j(x)=o(\norm x).
	\]
	For this fixed $a$, the number $q_j$ may be small but is strictly
	positive. Thus sufficiently far out on this end,
	\[
	dh_a(\partial_1)=q_j+\partial_1v_j>\frac{q_j}{2}.
	\]
	The graph metric satisfies
	$\norm{\partial_1}_g^2=1+\norm{\partial_1u_j}^2\le4$
	there. By Cauchy--Schwarz, $dh_a(\partial_1)=g(\nabla h_a,\partial_1)\leq |\nabla h_a||\partial_1|$, so
	$\norm{\nabla h_a}\ge q_j/4$ on that tail.
	Since there are only finitely many ends, the minimum of these
	positive bounds is positive after enlarging a compact core.
	
	The regular-value condition gives the Morse property by the
	preceding argument. By $|\nabla h_a|\ge c_a>0$, all critical points belong to a
	compact subset. Each critical point is isolated by nondegeneracy. Hence, there are finitely many critical points.
	
	We focus on one end and omit
	the subscript $j$, writing $h_a\circ F_j^{-1}(x)=qx_1+v(x)$ for $\norm x>R_0$.
	Take a smooth cutoff $\eta:[0,\infty)\to[0,1]$ that is $0$
	on $[0,1]$ and $1$ on $[2,\infty)$. For $R>R_0$ define
	\[
	\widetilde v(x)=
	\begin{cases}
		0,&\norm x\le R,\\
		\eta(\norm x/R)v(x),&\norm x>R.
	\end{cases}
	\]
	This is smooth, and it equals $v$ when $\norm x\ge2R$.
	In the transition annulus,
	\[
	\norm{D\widetilde v}
	\le\norm{Dv}+\frac{C}{R}\norm{v},
	\]
	where $C$ depends only on $\eta$. Both terms are uniformly
	small for large $R$, since $Dv\to0$ and $v=o(\norm x)$;
	the derivative is zero inside the inner disk and equals $Dv$
	outside the outer disk. Choose $R$ so large that
	$\sup_{\R^n}\norm{D\widetilde v}<q/2$.
	
	Define
	\[
	\widetilde h(x)=qx_1+\widetilde v(x),\qquad
	\Phi(x_1,x')=(\widetilde h(x_1,x'),x').
	\]
	For each fixed $x'$, the function
	$s\mapsto qs+\widetilde v(s,x')$ has derivative greater
	than $q/2$, and hence is strictly increasing with limits
	$\pm\infty$ as $s\to\pm\infty$. It is therefore a bijection
	of $\R$, with smooth inverse by the inverse function theorem.
	Thus $\Phi$ is a global diffeomorphism of $\R^n$, and
	\[
	\{\widetilde h\le t\}
	=\Phi^{-1}((-\infty,t]\times\R^{n-1})
	\cong(-\infty,t]\times\R^{n-1}
	\]
	for every $t\in\R$.
	
	Perform this construction on every end, possibly with different
	$R_j$, and extend the compact core $K$ to contain $\norm x\le2R_j$. The original height agrees with
	$\widetilde h_j$ on the remaining $j$th tail. Because $K$ and
	the finitely many disks $\overline B_{2R_j}$ are compact,
	we can choose $t$ below
	\[
	\min_K h_a,\qquad
	\min_{\overline B_{2R_j}}\widetilde h_j
	\quad(1\le j\le e).
	\]
	Then $M_a^t$ misses $K$, and each full extended sublevel
	$\{\widetilde h_j\le t\}$ misses its inner disk. Hence the
	$j$th tail contributes exactly one nonempty component,
	diffeomorphic to a closed half-space. There are precisely $e$
	such components.
	
	Similarly, choose $T$ above
	$\max_K h_a$ and above every
	$\max_{\overline B_{2R_j}}\widetilde h_j$.
	Then $K\subset M_a^T$, and every inner disk is contained in
	the corresponding extended sublevel. Given any point of an end
	sublevel, the connected set $\{\widetilde h_j\le T\}$ contains
	a path from that point to the center of its disk. Stop this
	path at its first intersection with the disk boundary.
	The preceding part lies entirely on the original end and
	satisfies $h_a\le T$; its endpoint lies on the boundary
	already absorbed into $K$. Thus every point of $M_a^T$
	can be connected within $M_a^T$ to the connected core $K$.
\end{proof}

\subsection{Relative homology of a low sublevel}\label{sec:relative-homology}
Write $m_k(f)$ for the number of critical points of index $k$ of a
Morse function $f$. We use homology with coefficients in $\Z_2$.

\begin{proposition}\label{prop:morse-applicability}
Let $f=h_a$ be a height as in Lemma~\ref{lem:generic-height}. For
every sufficiently negative $t$, put $L=M_a^t$. Then
\begin{equation}\label{eq:relative-morse-input}
 m_k(f)\ge \dim H_k(M,L;\Z_2)\qquad(0\le k\le n).
\end{equation}
\end{proposition}

\begin{proof}
	Write $M^s=\{x\in M:f(x)\le s\}$. By
	Lemma~\ref{lem:generic-height}, $f$ is Morse with finitely many
	critical points, and there exist a compact set $K\subset M$ and
	$c_a>0$ such that
	\[
	|\nabla f|\ge c_a\qquad\text{on }M\setminus K.
	\]
	Choose $t$ sufficiently negative that $L=M^t$ is nonempty and
	$t$ lies below every critical value of $f$.
	
	We first verify that the standard Morse deformation argument applies
	on the noncompact manifold $M$. Suppose $u<v$ and $f$ has no critical
	points in $f^{-1}([u,v])$. Then
	\[
	\inf_{f^{-1}([u,v])}|\nabla f|>0.
	\]
	Indeed, the gradient has the lower bound $c_a$ outside $K$, while
	$K\cap f^{-1}([u,v])$ is compact and contains no critical point.
	On this band consider the normalized negative gradient field
	$-\frac{\nabla f}{|\nabla f|^2}$.
	It has bounded norm, and along each of its integral curves
	$\Phi_\tau$ we have
	\[
	\frac{d}{d\tau}f(\Phi_\tau(x))=-1,
	\qquad f(\Phi_\tau(x))=f(x)-\tau.
	\]
	Construct strong deformation retraction on $M^v$ by
	\[
	H(s,x)=
	\begin{cases}
		x,&f(x)\le u,\\
		\Phi_{s(f(x)-u)}(x),&u<f(x)\leq v,
	\end{cases}
	\qquad0\le s\le1,
	\]
	which contracts $M^v$ onto $M^u$.
	
	Choose $T>\max_K f$ with $T>t$. Every critical point lies in $K$,
	so $T$ is above all critical values, and $|\nabla f|\ge c_a$ on
	$\{f\ge T\}$. The same normalized negative gradient flow, now run
	for time $f(x)-T$ when $f(x)>T$, gives a strong deformation
	retraction $M\searrow M^T$. It fixes $M^T$, and in particular $L$.
	Thus 
	\begin{equation}\label{eq:morse-upper-retraction}
		H_k(M^T,L;\Z_2)\cong H_k(M,L;\Z_2)
		\qquad(0\le k\le n).
	\end{equation}
	
	We next examine a critical level $c$. Choose regular values
	$u<c<v$ such that $c$ is the only critical value in $[u,v]$.
	Take pairwise disjoint Morse coordinate neighborhoods of the
	critical points at level $c$. On the complement of smaller such
	neighborhoods in $f^{-1}([u,v])$, the gradient is bounded away from
	zero: this follows from the compactness argument above, together
	with the estimate outside $K$. Hence, the normalized gradient
	flow controls the part of the band away from the critical points,
	and the usual local Morse handle-attachment construction applies
	(cf.\ \cite[Theorem~2.2]{vanGemmeren1998}). In a Morse chart at a
	critical point $x$ of index $\lambda=\ind_f(x)$,
	\[
	f(y)=c-y_1^2-\cdots-y_\lambda^2
	+y_{\lambda+1}^2+\cdots+y_n^2.
	\]
	Crossing $c$ attaches one $\lambda$-handle for each such point.
	Equivalently, up to homotopy, $M^v$ is obtained from $M^u$ by
	attaching one cell of dimension $\ind_f(x)$ for each critical
	point $x$ with $f(x)=c$. In particular,
	\begin{equation}\label{eq:morse-critical-band}
		\dim H_k(M^v,M^u;\Z_2)
		=\#\{x\in\operatorname{Crit}(f):f(x)=c,\ \ind_f(x)=k\}.
	\end{equation}
	
	Let $c_1<\cdots<c_s$ be the distinct critical values of $f$.
	If $s=0$, the no-critical-level deformation above gives
	$M^T\searrow L$, and the conclusion follows from
	\eqref{eq:morse-upper-retraction}. Otherwise choose regular values
	\[
	t=t_0<t_1<\cdots<t_s=T
	\]
	so that $(t_{j-1},t_j)$ contains precisely $c_j$. Set
	$A_j=M^{t_j}$. By \eqref{eq:morse-critical-band},
	\[
	\dim H_k(A_j,A_{j-1};\Z_2)
	=\#\{x\in\operatorname{Crit}(f):f(x)=c_j,\ \ind_f(x)=k\}.
	\]
	The exact sequence of the triple $L\subset A_{j-1}\subset A_j$
	contains
	\[
	H_k(A_{j-1},L;\Z_2)\longrightarrow
	H_k(A_j,L;\Z_2)\longrightarrow
	H_k(A_j,A_{j-1};\Z_2).
	\]
	Its exactness implies
	\[
	\dim H_k(A_j,L;\Z_2)
	\le \dim H_k(A_{j-1},L;\Z_2)
	+\dim H_k(A_j,A_{j-1};\Z_2).
	\]
	Since $A_0=L$, summing over $j=1,\ldots,s$ gives
	\[
	\begin{aligned}
		\dim H_k(M^T,L;\Z_2)
		&\le \sum_{j=1}^s
		\#\{x\in\operatorname{Crit}(f):
		f(x)=c_j,\ \ind_f(x)=k\}\\
		&=m_k(f).
	\end{aligned}
	\]
	Together with \eqref{eq:morse-upper-retraction}, this proves
	\eqref{eq:relative-morse-input} for every $0\le k\le n$.
\end{proof}

\begin{proposition}\label{prop:relative-morse}
For each direction in Lemma~\ref{lem:generic-height},
\begin{equation}\label{eq:relative-morse}
 m_1(h_a)\ge b_1(M)+e-1,\qquad
 m_k(h_a)\ge b_k(M)\quad(2\le k\le n-1).
\end{equation}
\end{proposition}

\begin{proof}
A sufficiently low $L$ has $e$ contractible components, so
$H_j(L)=0$ for $j>0$ and $H_0(L)=\Z_2^e$. The exact sequence of
the pair $(M,L)$ 
\[
 \cdots\to H_k(L)\to H_k(M) \to H_k(M,L)\to H_{k-1}(L)\to H_{k-1}(M)\to \cdots
\]
gives $H_k(M,L)\cong H_k(M)$ for $k\ge2$, and
\begin{equation}\label{eq:low-relative-exact}
 0\longrightarrow H_1(M)\longrightarrow H_1(M,L)
 \longrightarrow\Z_2^e\xrightarrow{\ \sum\ }\Z_2\longrightarrow0 .
\end{equation}
The last map sums the coordinates since all components of $L$ lie
in the same component of $M$. Its kernel has dimension $e-1$, hence
$\dim H_1(M,L)=b_1(M)+e-1$. Substitution into
\eqref{eq:relative-morse-input} gives the bounds.
\end{proof}

\subsection{Counting critical points with the Gauss map}\label{sec:gauss-count}
For $1\le k\le n-1$, let $U_k$ be the open subset of $N^1M$ where
$S_\nu$ is nonsingular of index $k$. The sets $U_k$ are disjoint
in $N^1M$. Set
\[
 I_k=\int_{U_k}|\det S_\nu|\dd V,\qquad
 \beta_1=b_1(M)+e-1,\qquad
 \beta_k=b_k(M)\quad(2\le k\le n-1).
\]

\begin{proposition}\label{prop:gauss-count}
Each $I_k$ is finite and
\begin{equation}\label{eq:ik-count}
 I_k=\int_{\Sn^{N-1}}m_k(h_a)\dd a\ge\Vol(\Sn^{N-1})\beta_k.
\end{equation}
In particular,
\begin{equation}\label{eq:index1-count}
 I_1\ge\Vol(\Sn^{N-1})(b_1(M)+e-1)\ge(e-1)\Vol(\Sn^{N-1}).
\end{equation}
\end{proposition}

\begin{proof}
	We first establish the index-counting identity. Recall that
	the normal Gauss map is
	\[
	G:N^1M\longrightarrow\Sn^{N-1},
	\qquad G(x,\nu)=\nu.
	\]
	By \eqref{eq:hessian-height}, a point $x\in M$ is critical
	for $h_a$ if and only if $a\in N_xM$. At such a point,
	the Hessian of $h_a$ is represented by the shape operator
	$S_a$. Consequently,
	\[
	G^{-1}(a)
	=\{(x,a):x\in\operatorname{Crit}(h_a)\},
	\]
	and a critical point $x$ is nondegenerate of Morse index
	$k$ precisely when $(x,a)\in U_k$.
	
	We recall why the Jacobian of $G$ is $|\det S_\nu|$.
	The normal connection gives orthogonal decompositions
	\[
	T_{(x,\nu)}N^1M
	\cong T_xM\oplus(N_xM\cap\nu^\perp)
	\]
	and
	\[
	T_\nu\Sn^{N-1}
	=T_xM\oplus(N_xM\cap\nu^\perp).
	\]
	For a horizontal tangent vector $X\in T_xM$, the
	Weingarten formula gives
	$dG(X,0)=-S_\nu X$, while for a vertical vector
	$\eta\in N_xM\cap\nu^\perp$, one has $dG(0,\eta)=\eta$.
	Thus
	\[
	dG_{(x,\nu)}
	=
	\begin{pmatrix}
		-S_\nu&0\\
		0&I_{p-1}
	\end{pmatrix}
	\]
	with respect to these decompositions. In particular,
	as in \eqref{eq:gauss-jacobian},
	\[
	J_G(x,\nu)=|\det dG_{(x,\nu)}|
	=|\det S_\nu|.
	\]
	Hence $G$ is a local diffeomorphism on $U_k$.
	
	Since $U_k$ is an open subset of the smooth manifold
	$N^1M$, choose an increasing compact exhaustion
	$K_j\subset U_k$, $j\ge1$. Applying the area formula
	\cite[\S3.2.3]{Federer1969} to the restriction of $G$
	to $K_j$ yields
	\[
	\int_{K_j}|\det S_\nu|\dd V
	=
	\int_{\Sn^{N-1}}
	\#\bigl(G^{-1}(a)\cap K_j\bigr)\dd a.
	\]
	Both integrands are nonnegative. As $j\to\infty$,
	monotone convergence therefore gives
	\[
	\int_{U_k}|\det S_\nu|\dd V
	=
	\int_{\Sn^{N-1}}
	\#\bigl(G^{-1}(a)\cap U_k\bigr)\dd a.
	\]
	
	By Sard's theorem, almost every $a\in\Sn^{N-1}$
	is a regular value of $G$. Moreover,
	Lemma~\ref{lem:generic-height} implies that,
	outside a set of measure zero, $h_a$ is a Morse
	function with finitely many critical points.
	For every direction satisfying both properties,
	the preceding correspondence gives
	\[
	\#\bigl(G^{-1}(a)\cap U_k\bigr)=m_k(h_a).
	\]
	We thus obtain the index-counting form of the
	Chern--Lashof formula
	\cite{ChernLashof,ShiohamaXu1997,ShiohamaXu2011}:
	\begin{equation}\label{eq:area-index}
		I_k
		=\int_{\Sn^{N-1}}
		\#\bigl(G^{-1}(a)\cap U_k\bigr)\dd a
		=\int_{\Sn^{N-1}}m_k(h_a)\dd a.
	\end{equation}
	
	We next show that $I_k$ is finite. Let
	$\lambda_1,\ldots,\lambda_n$ be the eigenvalues of
	$S_\nu$. By the arithmetic--geometric mean
	inequality,
	\[
	\begin{aligned}
		|\det S_\nu|
		&=\left(\prod_{i=1}^n\lambda_i^2\right)^{1/2}\\
		&\le
		\left(\frac1n\sum_{i=1}^n\lambda_i^2\right)^{n/2}\\
		&=n^{-n/2}|S_\nu|^n.
	\end{aligned}
	\]
	For an orthonormal basis $e_1,\ldots,e_n$ of $T_xM$,
	the definition of the shape operator and
	Cauchy--Schwarz give
	\[
	\begin{aligned}
		|S_\nu|^2
		&=\sum_{i,j=1}^n
		\langle S_\nu e_i,e_j\rangle^2\\
		&=\sum_{i,j=1}^n
		\langle A(e_i,e_j),\nu\rangle^2\\
		&\le\sum_{i,j=1}^n|A(e_i,e_j)|^2
		=|A|^2.
	\end{aligned}
	\]
	Thus, for every $(x,\nu)\in N^1M$,
	\[
	|\det S_\nu|\le n^{-n/2}|A(x)|^n.
	\]
	Using the normal-fiber integration formula
	\eqref{eq:normal-fiber-integral}, we obtain
	\[
	\begin{aligned}
		I_k
		&\le\int_{N^1M}|\det S_\nu|\dd V\\
		&\le n^{-n/2}
		\int_M\int_{S(N_xM)}
		|A(x)|^n\dd\nu\dd\mu(x)\\
		&=n^{-n/2}\Vol(\Sn^{p-1})
		\int_M|A|^n\dd\mu\\
		&<\infty.
	\end{aligned}
	\]
	For $p=1$, the integral over $\Sn^0$ is understood
	with counting measure, so the same argument applies.
	
	Finally, Proposition~\ref{prop:relative-morse}
	gives, for almost every $a\in\Sn^{N-1}$,
	\[
	m_k(h_a)\ge\beta_k,
	\qquad 1\le k\le n-1.
	\]
	Integrating this inequality and applying
	\eqref{eq:area-index}, we conclude that
	\[
	\begin{aligned}
		I_k
		&=\int_{\Sn^{N-1}}m_k(h_a)\dd a\\
		&\ge\int_{\Sn^{N-1}}\beta_k\dd a\\
		&=\Vol(\Sn^{N-1})\beta_k.
	\end{aligned}
	\]
	This proves \eqref{eq:ik-count}. In particular,
	since $\beta_1=b_1(M)+e-1$ and $b_1(M)\ge0$,
	\[
	\begin{aligned}
		I_1
		&\ge\Vol(\Sn^{N-1})(b_1(M)+e-1)\\
		&\ge(e-1)\Vol(\Sn^{N-1}),
	\end{aligned}
	\]
	which is \eqref{eq:index1-count}.
\end{proof}

\begin{corollary}\label{cor:paired-index}
For $1\le k<n/2$,
\begin{equation}\label{eq:paired-indices}
 I_k=I_{n-k}\ge\Vol(\Sn^{N-1})\max\{\beta_k,\beta_{n-k}\}.
\end{equation}
When $n$ is even,
\begin{equation}\label{eq:middle-index}
 I_{n/2}\ge\Vol(\Sn^{N-1})b_{n/2}(M).
\end{equation}
\end{corollary}

\begin{proof}
The map $(x,\nu)\mapsto(x,-\nu)$ preserves $\dd V$ and
$|\det S_\nu|$. Since $S_{-\nu}=-S_\nu$, this map exchanges $U_k$
and $U_{n-k}$. Applying \eqref{eq:ik-count} to both indices gives
\eqref{eq:paired-indices}. For $k=n/2$,  \eqref{eq:ik-count} gives \eqref{eq:middle-index} directly.
\end{proof}

\section{Curvature inequalities and proof of Theorem \ref{main1}}\label{sec:algebra}

We establish optimal inequalities of a trace-free symmetric operator and then apply them to the proof of Theorem \ref{main1}. 
\begin{lemma}[Determinant at fixed index]\label{lem:all-det}
Let $T$ be an invertible trace-free symmetric operator on $\R^n$
with $k$ negative eigenvalues. Then
\begin{equation}\label{eq:all-det}
 |T|^n\ge c_{n,k}|\det T|,\qquad
 c_{n,k}=n^{n/2}\left(\frac{n-k}{k}\right)^{n/2-k}.
\end{equation}
The optimal coefficient satisfies $c_{n,k}=c_{n,n-k}$. Equality
holds precisely when the eigenvalues of each sign are equal.
\end{lemma}
\begin{proof}
Write the eigenvalues as $-x_1,\ldots,-x_k,y_1,\ldots,y_\ell$,
where $\ell=n-k$ and $x_i,y_j>0$. Since $\tr T=0$, we suppose $\sum_i x_i=\sum_j y_j=s$.
By Cauchy--Schwarz and the arithmetic--geometric mean inequality,
\[
 |T|^2=\sum_i x_i^2+\sum_j y_j^2 \ge s^2\left(\frac1k+\frac1\ell\right)
       =\frac{ns^2}{k\ell},\qquad
 |\det T|=x_1\cdots x_k\cdot y_1\cdots y_l \le \left(\frac sk\right)^k
                \left(\frac s\ell\right)^\ell .
\]
Divide the first by the second inequality,
\[
 \frac{|T|^n}{|\det T|}
 \ge n^{n/2}k^{k-n/2}\ell^{\ell-n/2}
 = n^{n/2}\left(\frac\ell k\right)^{n/2-k}.
\]
Both inequalities are equalities precisely for $x_i=s/k$ and
$y_j=s/\ell$. Replacing
$T$ by $-T$ exchanges the two indices.
\end{proof}

\begin{corollary}[The index-one coefficient]\label{lem:index-one-det}
For $k=1$ or $n-1$, \eqref{eq:all-det} reads
\begin{equation}\label{eq:index-one-det}
 |T|^n\ge c_{n,1}|\det T|,\qquad
 c_{n,1}=c_n=\frac{[n(n-1)]^{n/2}}{n-1}.
\end{equation}
The equality spectra are the nonzero multiples of
$(-(n-1),1,\ldots,1)$, up to permutation.
\end{corollary}
\begin{proof}
Apply Lemma~\ref{lem:all-det} to $T$ or $-T$ with $k=1$.
The trace equation fixes the remaining eigenvalue.
\end{proof}

\begin{lemma}[Averaging over the normal sphere]\label{lem:normal-average}
Let $M^n$ be immersed in $\R^{n+p}$, with $n\ge3$ and $N=n+p$.
At each point $x$,
\begin{equation}\label{eq:normal-average}
 \int_{S(N_xM)}|S_\nu|^n\dd\nu
 \le \frac{2\Vol(\Sn^{N-1})}{\Vol(\Sn^n)}|A|^n .
\end{equation}
The coefficient is optimal. At $A(x)\ne0$, equality holds precisely
when $N_1(x)=\operatorname{span}\{A(X,Y):X,Y\in T_xM\}$ is a line.
\end{lemma}
\begin{proof}
For any point $x$, choose an orthonormal normal frame $\{\nu_1,\cdots,\nu_p\}$. For a fixed unit normal vector $\nu$, suppose 
\[
 \nu = \sum_{\alpha=1}^p t_\alpha \nu_\alpha,~\sum_{\alpha=1}^p t^2_\alpha=1,
\]
then $S_\nu=\sum t_\alpha S_{\nu_\alpha}$. Hence, 
\[
 |S_\nu|^2=\left\langle \sum_\alpha t_\alpha S_{\nu_\alpha},\sum_\beta t_\beta S_{\nu_\beta} \right\rangle=\sum_{\alpha,\beta} t_\alpha t_\beta\langle S_{\nu_\alpha},S_{\nu_\beta} \rangle=\sum_{\alpha,\beta} t_\alpha t_\beta \tr(S_{\nu_\alpha}S_{\nu_\beta}).
\]

Set $Q_{\alpha\beta}=\tr(S_{\nu_\alpha}S_{\nu_\beta})$ and the matrix $Q=(Q_{\alpha\beta})$, then 
\[
 |S_\nu|^2 = t^T Q t,
\]
hence, the matrix $Q$ is semi-definite positive. 
Diagonalize it by choosing the normal frame. Its eigenvalues
$\lambda_\alpha\ge0$ satisfy
\begin{equation}\label{eq:normal-average-gram}
 |S_\nu|^2=\sum_{\alpha=1}^p\lambda_\alpha t_\alpha^2,\qquad
 \sum_{\alpha=1}^p\lambda_\alpha=|A|^2,\qquad
 \nu=\sum_{\alpha=1}^p t_\alpha\nu_\alpha .
\end{equation}
At $A=0$ both sides of \eqref{eq:normal-average} vanish. Otherwise
put $w_\alpha=\lambda_\alpha/|A|^2$. Since $\sum w_\alpha=1$ and
$n/2>1$, by Jensen inequality,
\begin{equation}\label{eq:normal-average-jensen}
 |S_\nu|^n=|A|^n\left(\sum_\alpha w_\alpha t_\alpha^2\right)^{n/2}
 \le |A|^n\sum_\alpha w_\alpha|t_\alpha|^n.
\end{equation}
Integrate on both sides,
\begin{equation}\label{eq:normal-average-integ}
 \int_{S(N_xM)} |S_\nu|^n \le |A|^n\sum_\alpha w_\alpha\int_{\mathbb{S}^{p-1}} |t_\alpha|^n=|A|^n\int_{\mathbb{S}^{p-1}} |t_1|^n,
\end{equation}
since $\int_{\mathbb{S}^{p-1}} |t_\alpha|^n=\int_{\mathbb{S}^{p-1}} |t_1|^n$ for each $\alpha$ and $\sum w_\alpha=1$.

For $p\ge2$,
\begin{equation}\label{eq:normal-average-moment}
\begin{aligned}
 J_{n,p}
 &=\int_{\Sn^{p-1}}|t_1|^n\dd\nu
   =\Vol(\Sn^{p-2})\int_{-1}^1|t|^n(1-t^2)^{(p-3)/2}\dd t\\
 &=\Vol(\Sn^{p-2}) B\!\left(\frac{n+1}{2},\frac{p-1}{2}\right)
   =\frac{2\Vol(\Sn^{N-1})}{\Vol(\Sn^n)}.
\end{aligned}
\end{equation}
For $p=1$ the counting measure on $\Sn^0$ gives $J_{n,1}=2$,
this equal to the value of the left hand side. Combining
\eqref{eq:normal-average-integ} proves the inequality.

If there are two weights $w_\alpha,w_\beta$ being positive, strict convexity
makes \eqref{eq:normal-average-jensen} strict inequality on the open set
$t_\alpha^2\ne t_\beta^2$, which has positive spherical measure.
If only one weight $w_\alpha$ is positive, it gives equality everywhere. Thus equality
at a nonflat point is equivalent to $\rank Q=1$. 

Define $L:~N_xM\to \mathrm{Sym}(T_x M)$ by $L(\nu)=S_\nu$, then
\begin{eqnarray*}
	\nu\in\mathrm{ker}L &\Longleftrightarrow& S_\nu=0\\
	&\Longleftrightarrow& \langle S_\nu X,Y \rangle,~\forall X,Y\in T_xM\\
	&\Longleftrightarrow& \langle A(X,Y),\nu \rangle=0,~\forall X,Y\in T_xM.
\end{eqnarray*}
This shows that $\mathrm{ker}L=N_1(x)^\perp$.

For any normal vector $\nu=\sum t_\alpha \nu_\alpha$ and $t^T=(t_1,\cdots,t_\alpha)$, since $Q$ is semi-positive, 
\[
 t^TQ=0\Longleftrightarrow t^TQt=0 \Longleftrightarrow |S_\nu|^2=0 \Longleftrightarrow S_\nu=0.
\]
Hence $\mathrm{ker}Q=\ker L=N_1(x)^\perp$. Then equality holding in (\ref{eq:normal-average}) is equivalent to 
\[
 1=\rank Q=p-\dim \ker Q=p-\dim N_1(x)^\perp=\dim N_1(x).
\]
\end{proof}

\subsection{Proof of Theorem \ref{main1}}\label{sec:full-betti}

\begin{proof}[Proof of Theorem \ref{main1}]
Minimality gives $\tr S_\nu=\ip{\tr_g A}{\nu}=0$, so each nonsingular
$S_\nu$ belongs to exactly one of $U_1,\ldots,U_{n-1}$.
Integrate \eqref{eq:all-det} on each $U_k$ by replacing $T=S_\nu$, then use
\eqref{eq:normal-average} and the fiber integration formula
\eqref{eq:normal-fiber-integral}. This gives
\begin{eqnarray}\label{eq:full-chain}
 \frac{2\Vol(\Sn^{N-1})}{\Vol(\Sn^n)}\int_M |A|^n\dd\mu
 &\ge&\int_{N^1M}|S_\nu|^n\dd V\\\nonumber
 &\ge&\sum_{k=1}^{n-1}\int_{U_k} |S_\nu|^n\dd V\\\nonumber
 &\ge&\sum_{k=1}^{n-1}c_{n,k}\int_{U_k} |\det S_\nu|\dd V\\\nonumber
 &=&\sum_{k=1}^{n-1}c_{n,k}I_k.\nonumber
\end{eqnarray}

For $k<n/2$, antipodal symmetry and
Corollary~\ref{cor:paired-index} gives
\[
 c_{n,k}I_k+c_{n,n-k}I_{n-k}
 =2c_{n,k}I_k
 \ge2c_{n,k}\Vol(\Sn^{N-1})\max\{\beta_k,\beta_{n-k}\}.
\]
When $n$ is even there will be the middle index, and its contribution
is at least $n^{n/2}\Vol(\Sn^{N-1})b_{n/2}$.
With $b_k=b_k(M)$ and $e=e(M)$, dividing the above by $\frac{2\Vol(\Sn^{N-1})}{\Vol(\Sn^n)}$ gives
\begin{equation}\label{eq:symmetric-betti}
\begin{aligned}
 \int_M|A|^n\dd\mu
 &\ge\gamma_n\max\{b_1+e-1,b_{n-1}\}\\
 &\quad+\Vol(\Sn^n)\biggl[
 \sum_{2\le k<n/2} c_{n,k}\max\{b_k,b_{n-k}\}
 +\begin{cases}
 \frac12 n^{n/2}b_{n/2}, &n\text{ even},\\
 0,&n\text{ odd}.
 \end{cases}\biggr]
\end{aligned}
\end{equation}
Here $c_{n,1}\Vol(\Sn^n)=\gamma_n$.
By Lemma~\ref{rem:betti-symmetry}, $b_{n-1}=b_1+e-1$ and
$b_{n-k}=b_k$ for $2\le k\le n-2$. Substitution gives
\eqref{eq:main-topology}. 
The catenoid has $e=2$ and $b_1=0$ and has energy $\gamma_n$
by \eqref{eq:catenoid-energy}. This proves sharpness for each $p$.
\end{proof}

For the equality case of the index-one bound, retain just the terms
at indices $1$ and $n-1$ in \eqref{eq:full-chain}. This gives
\begin{equation}\label{eq:indexone-chain}
 \frac{2\Vol(\Sn^{n+p-1})}{\Vol(\Sn^n)}\int_M |A|^n\dd\mu
 \ge\int_{N^1M}|S_\nu|^n\dd V
 \ge 2c_{n,1}I_1
 \ge \frac{2\Vol(\Sn^{n+p-1})}{\Vol(\Sn^n)}\gamma_n\beta_1 .
\end{equation}
Since $\frac{2\Vol(\Sn^{n+p-1})}{\Vol(\Sn^n)}\gamma_n=2c_{n,1}\Vol(\Sn^{N-1})$, equality in
\eqref{eq:topology-ends} forces equality in each step of
\eqref{eq:indexone-chain}, in particular in the normal-sphere integral.

\subsection{Comparison with the unweighted estimate}\label{sec:unweighted}

Let $B(M)=\sum_{k=0}^n b_k(M)$ and let $\tau(i)$ denote the normalized
total absolute curvature. Since $\det S_\nu$ vanishes off the union of
the $U_k$, the count in \eqref{eq:ik-count} gives
\begin{equation}\label{eq:planar-tau}
\begin{aligned}
 \tau(i)&=\frac{1}{\Vol(\Sn^{N-1})}\int_{N^1M}|\det S_\nu|\dd V
         =\frac{1}{\Vol(\Sn^{N-1})}\sum_{k=1}^{n-1}I_k\\
 &\ge b_1+e-1+\sum_{k=2}^{n-1}b_k=B(M)+e-2 .
\end{aligned}
\end{equation}
The last equality uses $b_0=1$ and $b_n=0$. For comparison, van
Gemmeren's general estimate for semi-Wintgen immersions with
$H$-spherical ends \cite[Theorem~3.3]{vanGemmeren1998} is
\begin{equation}\label{eq:vg98-tau}
 \tau(i)\ge
 \begin{cases}
 B(M)-1,&e=1,\\
 B(M)-2,&e\ge2 .
 \end{cases}
\end{equation}
The graphical ends have spherical cross sections, and their generic
heights satisfy the Morse--Palais--Smale condition of that theorem.
For $e\ge2$, the identification of one contractible low sublevel on each
end increases the lower bound for $\tau(i)$ by $e$.

The arithmetic--geometric mean inequality on the squares of the
eigenvalues gives $|S_\nu|^n\ge n^{n/2}|\det S_\nu|$ in every
signature. Integrating and applying \eqref{eq:normal-average} and
\eqref{eq:planar-tau} gives
\begin{equation}\label{eq:unweighted}
 \int_M|A|^n\dd\mu
 \ge\frac{n^{n/2}\Vol(\Sn^n)}{2}
 \left(b_1+e-1+\sum_{k=2}^{n-1}b_k\right).
\end{equation}

\section{Sharpness and rigidity}\label{sec:rigidity}

We compute the energy of the ordinary catenoid and its ends to show the sharpness in Theorem \ref{main1}. 
\subsection{The rotational catenoid}\label{sec:catenoid-model}

For sake of completeness, we recall the standard rotational $n$-catenoid, a classical
minimal hypersurface; see, for example, the two-end
hypersurface analysis of Schoen~\cite{Schoen1983}, or Tam-Zhou \cite[\S2]{TamZhou2009}.

Our calculation below fixes the scale-invariant
$\int\norm A^n$ normalization and verifies attainment
of the index-one coefficient in Theorem~\ref{main1}.
The model lies in an affine $\R^{n+1}\subset\R^{n+p}$. Adding zero coordinates changes neither its induced metric nor the norm of its second fundamental form. Write the rotational immersion in meridian arclength as
\[
i(s,\theta)=(r(s)\theta,z(s)),
\qquad \theta\in\Sn^{n-1},\quad r>0,\quad (r')^2+(z')^2=1.
\]
Its induced metric is $ds^2+r(s)^2g_{\Sn^{n-1}}$, and a unit normal is $\nu=(-z'\theta,r')$. Differentiating in meridian gives principal curvatures
\[
\kappa_s=\langle (r''\theta,z''),(-z'\theta,r')\rangle=-r''z'+z''r'.
\]
For sphere directions, fixing $s$ gives a radius $r(s)$ sphere. Suppose $v\in T_\theta\mathbb{S}^{n-1}$ is a unit tangent vector, then the corresponding unit direction is $(v,0)$, along this direction, the second derivative is $-\frac{1}{r}(\theta,0)$, then  
\[
\kappa_\theta=\left\langle -\frac{1}{r}(\theta,0),\nu\right\rangle=\frac{z'}r,
\]
which is of multiplicity $n-1$.

For $r'\ne0$, differentiating the arclength identity gives $\kappa_s=z''/r'$. Then the minimal submanifold equation $H=0$ becomes
\[
z''+(n-1)\frac{r'}r z'=0,
\qquad\text{or equivalently}\qquad(r^{n-1}z')'=0.
\]

For $r'=0$, set $\varphi(s)$ such that $r'=\cos\varphi,~z'=\sin\varphi$, then $r''=-\varphi'\sin\varphi,~z''=\varphi'\cos\varphi$. Hence,
\[
 \kappa_s=-r''z'+z''r'=\varphi\sin^2\varphi+\varphi'\cos^2\varphi=\varphi',~\kappa_\theta=\frac{\sin\varphi}{r}.
\]
$H=0$ becomes $\varphi'+(n-1)\frac{\sin\varphi}{r}=0$, integrate $\cos\varphi=r'$ to get the same equation
\[
z''+(n-1)\frac{r'}r z'=0.
\]
Normalize the neck radius to one, put the neck at $s=0$, and choose the increasing $z$ direction. At the neck $r'=0$ and $z'=1$. Thus
\begin{equation}\label{eq:cat-first-integral}
	z'=r^{-(n-1)},\qquad
	(r')^2=1-r^{-2(n-1)},\qquad
	r''=(n-1)r^{1-2n}.
\end{equation}
The last equation with $r(0)=1$, $r'(0)=0$ gives a smooth even function $r(s)$, since $r(-s)$ also satisfies the initial valued equation (\ref{eq:cat-first-integral}). Since $r''=(n-1)r^{1-2n}>0$ and $r'(0)=0$, $r'>0$ increases strictly for $s>0$ and $|r'|<1,~r'(s)\to 1$ as $s\to+\infty$. Then $r$ tends to infinity as $s\to\infty$. Defining $z(0)=0$ by the first equation of (\ref{eq:cat-first-integral}) gives an odd function $z$. Hence the immersion is smooth across the neck and complete on $\R\times\Sn^{n-1}$. Since
\[
 \frac{dz}{dr}=\frac{dz/ds}{dr/ds}=\frac{r^{-(n-1)}}{\sqrt{1-r^{-2(n-1)}}}=\frac{1}{\sqrt{r^{2(n-1)}-1}},
\]
in graphical coordinates its two ends are
\begin{equation}\label{eq:cat-graphs}
	i_\pm(r,\theta)=\left(r\theta,\ \pm\int_1^r\frac{dt}{\sqrt{t^{2(n-1)}-1}}\right),
	\qquad r\ge1,
\end{equation}

Equations~\eqref{eq:cat-first-integral} and minimality give
\begin{equation}\label{eq:cat-spectrum}
	\kappa_\theta=r^{-n},\qquad \kappa_s=-(n-1)\kappa_\theta=-(n-1)r^{-n}.
\end{equation}
In particular $\norm A^2=n(n-1)r^{-2n}$. On either half,
\[
d\mu=r^{n-1}\,ds\dd\theta
=\frac{r^{n-1}}{\sqrt{1-r^{-2(n-1)}}}\dd r\dd\theta.
\]
Substituting $t=r^{-(n-1)}$, one obtains
\begin{equation}\label{eq:catenoid-energy}
	\begin{aligned}
		\int_{\Cat_n}\norm A^n\dd\mu
		&=2\int_{\Cat_n\cap\{s>0\}}\norm A^n\dd\mu\\
		&=2\Vol(\Sn^{n-1})[n(n-1)]^{n/2}
		\int_1^\infty\frac{r^{-n(n-1)-1}}{\sqrt{1-r^{-2(n-1)}}}\dd r\\
		&=\frac{2\Vol(\Sn^{n-1})[n(n-1)]^{n/2}}{n-1}
		\int_0^1\frac{t^{n-1}}{\sqrt{1-t^2}}\dd t\\
		&=\frac{[n(n-1)]^{n/2}}{n-1}
		\Vol(\Sn^{n-1})B\left(\frac n2,\frac12\right)=\frac{[n(n-1)]^{n/2}}{n-1}\Vol(\Sn^n).
	\end{aligned}
\end{equation}
Here $B(p,q)=\int_0^1v^{p-1}(1-v)^{q-1}\dd v$, and the last identity follows from $\Vol(\Sn^{n-1})B(n/2,1/2)=\Vol(\Sn^n)$, using the Gamma-function formula for sphere volumes. Scale invariance shows that every neck radius has the same energy.

By \eqref{eq:cat-first-integral}, $r'$ increases strictly from $-1$ to $1$. The normal
\[
\nu(s,\theta)=\bigl(-r(s)^{-(n-1)}\theta,r'(s)\bigr)
\]
therefore maps the catenoid diffeomorphically onto $\Sn^n$ minus its two poles: the last coordinate determines $s$, and the remaining coordinates determine $\theta$. Here the Gauss map is taken within its containing $\R^{n+1}$. Its topology is $\Sn^{n-1}\times\R$, so $e=2$, $b_1=0$, and $b_{n-1}=1$ for $n\ge3$, with all other positive-degree Betti numbers zero. It realizes equality in Theorem~\ref{main1}.

\subsection{Reduction of codimension and the proof of Theorem \ref{main2}}\label{sec:codim-reduction}

The local reduction below is the rank-one case of the classical
codimension-reduction principle \cite{Erbacher1971}. Here Codazzi
gives parallelism of the normal line, and harmonicity extends the
result across zeros of $A$.

\begin{proposition}[Reduction from a rank-one first normal space]\label{lem:codim-reduction}
Let $i:M^n\to\R^{n+p}$ be a connected minimal immersion.
Suppose the first normal space
$N_1(x)=\operatorname{span}\{A_x(X,Y):X,Y\in T_xM\}$
has dimension one on a nonempty open set. Then $i(M)$ lies
in an affine $(n+1)$-plane.
\end{proposition}

\begin{proof}
Choose a smooth unit section $\xi$ of the first normal line on a
connected coordinate neighborhood $V$ in this open set. Write
$A=h\xi$. Then $h$ is nonzero and trace free, so it has rank at
least two. 
For the normal connection $\nabla^\perp$, we calculate 
\begin{eqnarray*}
 (\nabla_X^\perp A)(Y,Z)&=&(\nabla_X^\perp )(A(Y,Z))-A(\nabla_XY,Z)-A(Y,\nabla_XZ)\\
 &=&X(h(Y,Z))\xi + h(Y,Z)\nabla_X^\perp\xi-h(\nabla_X Y,Z)\xi-h(Y,\nabla_XZ)\xi\\
 &=&(\nabla_X h)(Y,Z)\xi+h(Y,Z)\nabla_X^\perp\xi.
\end{eqnarray*}
Combing the Codazzi equation
\[
(\nabla_X^\perp A)(Y,Z)=(\nabla_Y^\perp A)(X,Z),
\]
we have 
\[
 (\nabla_X h)(Y,Z)\xi+h(Y,Z)\nabla_X^\perp\xi=(\nabla_Y h)(X,Z)\xi+h(X,Z)\nabla_Y^\perp\xi.
\]
Since $|\xi|=1$, $\nabla_X^\perp \xi\perp \xi$. Hence, the above equality project to the orthogonal
complement of $\xi$ gives
\begin{equation}\label{eq:codazzi-normal}
	h(Y,Z)\nabla_X^\perp\xi=h(X,Z)\nabla_Y^\perp\xi .
\end{equation}
For any point $x\in V$, if $h_x=0$, then $A_x(X,Y)=0$ which contradicts with $\dim N_1(x)=1$. Hence, $h\not=0$ on $V$. Since $h$ is trace-free, it has at least two non-zero eigenvalue. Choose an orthonormal frame $\{e_1,\cdots,e_n\}$ of $T_xM$ such that $h(e_i,e_j)=\lambda_i\delta_{ij}$, for any $i$, there is always some $j\not=i$ such that $\lambda_j\not=0$. Pick $X=e_i,~Y=e_j,~Z=e_j$ in (\ref{eq:codazzi-normal}), we get $\lambda_j\nabla_{e_i}^\perp \xi=0$. Thus,
\begin{equation}\label{eq:codim-parallel-line}
 \nabla^\perp\xi=0\quad\text{on }V.
\end{equation}

For any $x\in V$, define $E_x=di_x(T_x M)\oplus \mathbb{R}\xi(x)\subset\R^{n+p}$. Next, we want to prove there exists a fixed linear space $E$, such that $E_x=E,~\forall x\in V$.

The Gauss formula
$D_X Y=\nabla_X Y+h(X,Y)\xi$ shows $D_X Y(x)\in E_x$. The Weingarten formula $D_X\xi(x)=-S_\xi X(x)\in E_x$. For a local orthonormal frame over $V$, $\{u_1,\cdots,u_n\}$, and choose $u_{n+1}=\xi$, they form an orthonormal basis of $E_x$ for each $x\in V$. The above shows $D_X u_i(x)\in E_x$ for $i=1,\cdots,n+1$. Fix a point $x_0\in V$, for any point $x\in V$, connect them by a curve $\gamma(t)$ such that $\gamma(0)=x_0$ and $\gamma(1)=x$. Along this curve, set $U(t)=(u_1(t),\cdots,u_{n+1}(t))$, there exists a matrix $C(t)$ such that $U'(t)=U(t)C(t)$. Its solution is of the form $U(t)=U(0)P(t)$ for some invertible $P(t)$. Hence,
\[
 E_x=E_{\gamma(t)}=\mathrm{span}\{u_1(t),\cdots,u_{n+1}(t)\}=\mathrm{span}\{u_1(0),\cdots,u_{n+1}(0)\}=E_{\gamma(0)}=E_{x_0}:=E.
\]
Moreover, 
\[
 i(x)-i(x_0)=\int_0^1 di_{\gamma(t)}(\gamma'(t))\dd t\in E.
\]
Thus $i(V)\subset i(x_0)+E$, an affine $(n+1)$-plane.

Next, we show that $i(M)\subset i(x_0)+E$. For any fixed
$v\in E^\perp$, define $f_v(x)=\langle i(x)-i(x_0),v\rangle$
for $x\in M$. We have $f_v|_V=0$ and
$\Delta_M f_v(x)=\langle\Delta_M i(x),v\rangle=0$,
since $i$ is minimal. Thus $f_v$ is harmonic on $M$.
By unique continuation for harmonic functions on the connected
manifold $M$, $f_v\equiv0$. Hence $i(x)-i(x_0)\in E$
for all $x\in M$, i.e. $i(M)\subset i(x_0)+E$.
\end{proof}

\begin{proof}[Proof of Theorem~\ref{main2}]
	Suppose $\int_M |A|^n\dd \mu<\gamma_n$. 
    Furthermore, if $A\equiv0$, then $i(M)$ is an affine $n$-plane.
	
	Suppose now that $A\not\equiv0$.
	By Lemma~\ref{lem:one-end}, a nonflat
	complete minimal immersion of finite total
	curvature has at least two ends. Hence
	$e(M)\ge2$. Applying Theorem~\ref{main1},
	we obtain
	\begin{equation}\label{eq:thmB-lower-bound}
		\begin{aligned}
			\int_M |A|^n\dd \mu
			&\ge\gamma_n\bigl(b_1(M)+e(M)-1\bigr)\\
			&\ge\gamma_n,
		\end{aligned}
	\end{equation}
	which contradicts with the assumption that $\int_M |A|^n\dd \mu<\gamma_n$. 
	Hence, $i(M)$ is an affine $n$-plane.
	
	It remains to consider the equality case
	$\int_M |A|^n\dd \mu=\gamma_n$.
	Equality in \eqref{eq:thmB-lower-bound}
	implies
	\[
	b_1(M)+e(M)-1=1.
	\]
	Since $b_1(M)\ge0$ and $e(M)\ge2$, we have
	\begin{equation}\label{eq:thmB-topology-equality}
		e(M)=2,\qquad b_1(M)=0.
	\end{equation}
	
	We next show that the immersion reduces
	to codimension one.
	
	Recall the index-one inequality chain
	\eqref{eq:indexone-chain}:
	\[
	\frac{2\Vol(\Sn^{n+p-1})}{\Vol(\Sn^n)}\int_M |A|^n\dd \mu
	\ge
	\int_{N^1M}|S_\nu|^n\dd V
	\ge
	2c_{n,1}I_1
	\ge
	\frac{2\Vol(\Sn^{n+p-1})}{\Vol(\Sn^n)}\gamma_n
	\bigl(b_1(M)+e(M)-1\bigr),
	\]

	Since
	\[
	\int_M |A|^n\dd \mu
	=
	\gamma_n\bigl(b_1(M)+e(M)-1\bigr),
	\]
	the first and last terms are equal.
	Therefore equality holds throughout
	the chain. In particular,
	\[
	\int_{N^1M}|S_\nu|^n\dd V
	=
	\frac{2\Vol(\Sn^{n+p-1})}{\Vol(\Sn^n)}\int_M|A|^n\dd\mu.
	\]
	
	By the normal-fiber integration formula,
	this gives
	\[
	\int_M
	\left(
	\frac{2\Vol(\Sn^{n+p-1})}{\Vol(\Sn^n)}|A(x)|^n
	-
	\int_{S(N_xM)}|S_\nu|^n\dd\nu
	\right)\dd\mu(x)=0.
	\]
	By Lemma~\ref{lem:normal-average},
	\begin{equation}\label{eq:normal-averaging-defect}
	\frac{2\Vol(\Sn^{n+p-1})}{\Vol(\Sn^n)}|A(x)|^n
		-
		\int_{S(N_xM)}|S_\nu|^n\dd\nu\ge 0.
	\end{equation}
    Hence,
	\[
	\frac{2\Vol(\Sn^{n+p-1})}{\Vol(\Sn^n)}|A(x)|^n
	=\int_{S(N_xM)}|S_\nu|^n\dd\nu.
	\]
	
	Since $A\not\equiv0$, the set
	\[
	U=\{x\in M:|A(x)|>0\}
	\]
	is a nonempty open subset of $M$.
	By the equality characterization in
	Lemma~\ref{lem:normal-average}, we obtain
	\[
	N_1(x)=\operatorname{span}
	\{A_x(X,Y):X,Y\in T_xM\},\qquad
	\dim N_1(x)=1,\quad x\in U.
	\]

	Proposition~\ref{lem:codim-reduction}
	therefore implies that $i(M)$ lies
	in an affine $(n+1)$-plane
	$V\subset\R^{n+p}$.
	
	Identifying $V$ with $\R^{n+1}$,
	we may regard $i$ as a complete connected
	minimal hypersurface immersion with
	exactly two ends.
	By Proposition~\ref{prop:hypersurface-regularity},
	the immersion is regular at infinity.
	Schoen's classification theorem
	(Theorem~\ref{thm:schoen}) then implies
	that $i(M)$ is either an ordinary
	$n$-catenoid or a pair of affine
	hyperplanes. Since $M$ is connected, $i(M)$ must be an ordinary
	$n$-catenoid.
	
	Conversely, the computation
	\eqref{eq:catenoid-energy} shows that
	every ordinary $n$-catenoid satisfies
	\[
	\int_{\Cat_n}|A|^n\dd\mu
	=
	\frac{[n(n-1)]^{n/2}}{n-1}
	\Vol(\Sn^n)
	=
	\gamma_n.
	\]
	This proves the equality
	characterization and the sharpness.
\end{proof}
\section*{Acknowledgments}
The author thank Prof. Xingwang Xu for introducing the gap type theorem several years ago. Special thanks to Prof. Hongyu Wang and Prof. Peng Lu for their support.

\section*{Disclosure on AI assistance}
The principal mathematical ideas were formulated by the author. The original form of this paper is for codimension one, AI-assisted tools helped to find reference \cite{Erbacher1971} to formulate Proposition \ref{lem:codim-reduction} which made us be able to deal with the high codimensional case. Also they were used to help check selected computations of the manuscript. The author verified and completed all mathematical arguments and takes full responsibility for the content.

\end{document}